\documentclass[conference]{IEEEtran}
\IEEEoverridecommandlockouts
\usepackage{cite}
\usepackage{amsmath,amssymb,amsfonts}
\usepackage{algorithmic}
\usepackage{graphicx}
\usepackage{textcomp}
\usepackage{xcolor}
\usepackage{url}
\usepackage[numbers]{natbib}
\usepackage{hyperref}
\usepackage{comment}
\usepackage{placeins}
\usepackage{float}
\usepackage{enumerate}

\def\BibTeX{{\rm B\kern-.05em{\sc i\kern-.025em b}\kern-.08em
    T\kern-.1667em\lower.7ex\hbox{E}\kern-.125emX}}
\begin{document}

\title{Statistical data analysis for Tourism in Poland 
in R Programming Environment\\
}

\author{\IEEEauthorblockN{
Saad Ahmed Jamal}
\IEEEauthorblockA{\textit{MED-IIFA} \\
Universidade de Evora, Portugal \\
saad.jamal@uevora.pt}
}

\maketitle


\begin{abstract}
This study utilises the R programming language for statistical data analysis 
to understand Tourism dynamics in Poland. 
It focuses on methods for data visualisation, multivariate statistics, and hypothesis testing. To investigate the expenditure behavior of tourist, spending patterns, correlations, and associations among variables were analysed in the dataset. 
The results revealed a significant relationship between accommodation type and the purpose of trip,  showing that the purpose of a trip impacts the selection of accommodation.
A strong correlation was observed between organizer expenditure and private expenditure, indicating that individual spending are more when the spending on organizing the trip are higher.
However, no significant difference was observed in total expenditure across different accommodation types and purpose of the trip  
revealing that travelers tend to spend similar amounts regardless of their reason for travel or choice of accommodation.  
Although significant relationships were observed among certain variables, ANOVA could not be applied because the dataset was not able to hold on the normality assumption. In future, 
the dataset can be explored further to find more meaningful insights.
The developed code is available on GitHub: \url{https://github.com/SaadAhmedJamal/DataAnalysis_RProgEnv}.
\end{abstract}


\section{Introduction}

R is a high-level programming language widely used in statistical computing, data analysis, and geospatial applications. It was developed by \citet{0_ihaka1996r} at the University of Auckland and has since become a leading tool for statistical modeling and visualisation. The latest stable version, R ‘4.4.1,’ used in this study, is freely available at \url{https://cran.r-project.org/}. Its extensive package ecosystem allows researchers to perform efficient data wrangling, hypothesis testing, and visualization \citep{wickham2016getting, 4_wickham2017package}.  
There was a debate if R is any better and whether is showed be used for teaching purposes  \citep{1_gomes2018teaching}. A comprehensive learning tutorials were developed which made learning in R easier \citep{1b_de2018teaching}. With time, R has evolved be as important as Python in the field of Data Science.
R for statistical computing and data visualisation due to its flexibility and robust package ecosystem \citep{0_ihaka1996r}. The ability to perform exploratory data analysis (EDA), hypothesis testing, and regression modeling makes R an essential tool for empirical research \citep{wickham2016getting}.

Statistical data analysis plays a crucial role in tourism research and ecological studies, where identifying patterns, trends, and relationships in datasets can drive better decision-making \citep{ding2018analysis, davis2023writing}. For instance, in tourism analytics, understanding expenditure patterns, accommodation choices, and trip purposes enables businesses and policymakers to optimize tourism management strategies. Similarly, in ecological research, morphometric analyses—such as assessing the relationship between forearm length and body weight in bats—help scientists study species adaptation and environmental impacts \citep{5_harrell2019package}. Further these assessment can be used in regression and predictive analysis 
\citep{nusrat2024multiclassdepressiondetection}.

Tourism data analysis has been widely applied in economic studies to assess travel expenditure, seasonal trends, and consumer behavior \citep{ding2018analysis}. Statistical models 
provide insights into the factors influencing tourist spending patterns, allowing policymakers and businesses to optimize tourism planning and resource allocation.
\citet{ding2018analysis} examined domestic tourism consumption in China and found that economic indicators such as GDP and disposable income significantly influenced tourism expenditure. Statistical models such as ANOVA and t-tests have been widely used to assess variations among different categories of travelers, such as business versus leisure tourists \citep{milenkovski2019statistical}. Additionally, multivariate statistical techniques, including time-series models and clustering algorithms, have been employed to predict seasonal fluctuations in travel demand, assisting policymakers in optimizing tourism strategies.

Ecological research has extensively utilized statistical methods to analyse morphometric traits and their ecological significance. Studies indicate that biometric measurements, such as forearm length and body weight in bats, serve as key indicators of species adaptation and health \citep{davis2023writing}. Non-parametric statistical methods 
are particularly useful when normality assumptions are violated. \citet{davis2023writing} conducted a comparative analysis of statistical methodologies in ecological research, emphasising the importance of selecting appropriate methods based on data distribution. Furthermore, geometric morphometric techniques implemented in R have facilitated the quantification of shape variations and evolutionary patterns in ecological studies \citep{adams2013geomorph}. The existing literature highlights the effectiveness of statistical techniques in tourism and ecological research.

The primary goal of this study is to analyse the Tourism Dataset, which contains information on Polish household trips. The analysis aims to explore relationships between key variables, expenditure trends, and behavioral patterns in trip planning. Specifically, this study investigates:
\begin{enumerate}[i]
    \item The relationship between trip purpose, accommodation type, and nights spent on total expenditure.
    \item Hypothesis testing methodologies, including t-tests and ANOVA, to assess differences between groups.
    \item Graphical visualisations to represent data distributions and spending trends.
\end{enumerate}

Using hypothesis testing and data visualisation techniques, this study aims to contribute to the growing body of research on data-driven decision-making in tourism and consumer behavior analysis.
The integration of R-based data visualisation, hypothesis testing, and multivariate statistical methods has enhanced the reproducibility and accuracy of findings in the field. This study builds upon 
applying inferential statistical methods to analyse tourism expenditure and ecological morphometrics, contributing to the growing field of data-driven decision-making.

\section{Methodology}

\begin{figure*}[htbp]
\centerline{\includegraphics[width=\textwidth]{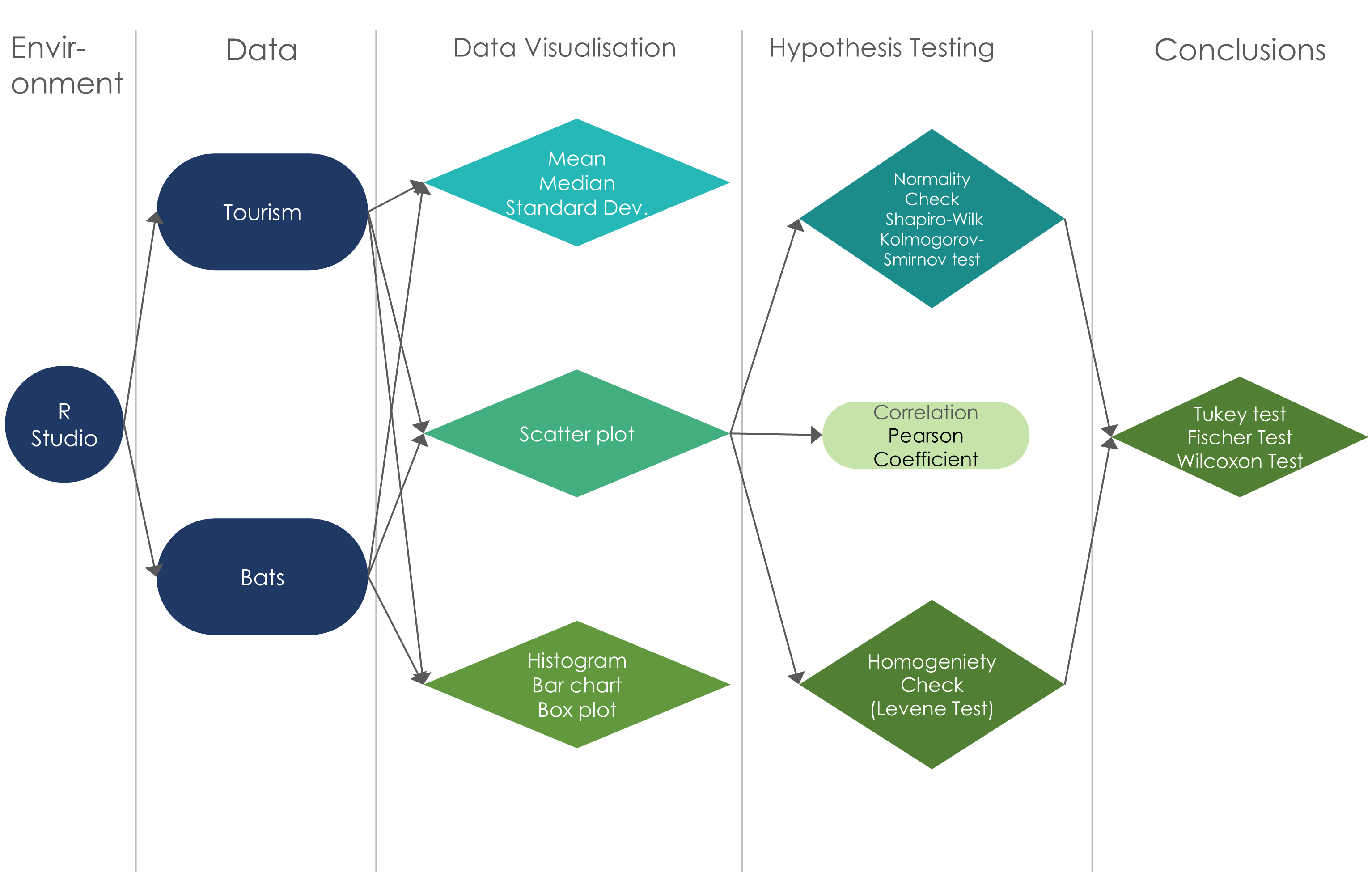}}
\caption{Methodological Flowchart.}
\label{figMethod}
\end{figure*}

A systematic, step-by-step approach was used for the analysis in R Studio.
Figure \ref{figMethod} shows the methodological flowchart for statistical modeling. 

\subsection{Packages Used}

Several R packages were utilised in this study for data handling, visualisation, and statistical analysis:

\begin{itemize}
    \item \texttt{readxl}: Used for reading Excel files (.xls and .xlsx) into R, allowing direct import of spreadsheet data without requiring external dependencies \citep{2_wickham2019package}. Unlike \texttt{read.csv()}, which handles CSV files, \texttt{readxl} supports native Excel formats.
    \item \texttt{dplyr}: A package for efficient data manipulation and transformation, part of the Tidyverse collection. It provides functions for filtering, summarizing, and modifying data frames \citep{3_yarberry2021dplyr}.
    \item \texttt{tidyverse}: A collection of R packages designed for data science, including \texttt{ggplot2}, \texttt{dplyr}, \texttt{tidyr}, \texttt{readr}, \texttt{tibble}, \texttt{stringr}, and \texttt{purrr} \citep{4_wickham2017package}. This suite simplifies data wrangling and visualisation.
    \item \texttt{ggplot2}: The most widely used package for data visualisation in R. It is based on the grammar of graphics, enabling the creation of layered plots \citep{wickham2016getting}.
    \item \texttt{reshape2}: Provides tools to reshape datasets between wide and long formats, making them more suitable for different types of statistical analyses \citep{wickham2007reshaping}.
    \item \texttt{Hmisc}: A package for statistical modeling and data analysis, offering functions for descriptive statistics, imputation, and regression modeling \citep{5_harrell2019package}.
\end{itemize}




\subsection{Dataset Description}


The dataset used in this study was collected from a sample survey of Polish households over a specific period. It is an open-source data available online \citep{Cellini2013-ul}. It consists of numerical and categorical variables, making it suitable for both descriptive and inferential statistical analyses. Table \ref{tab0} provides an overview of the dataset variables.

\begin{table}[H]
\caption{Dataset Variables and  Description}
\begin{center}
\begin{tabular}{|p{2.5cm}|p{6cm}|}
\hline 
\textbf{Variable} & \textbf{Description} \\
\hline
Year & Year in which the trip occurred. \\
PT (Purpose of Trip) & The reason for travel, such as leisure, business, or other. \\
AT (Accommodation Type) & Type of accommodation used during the trip (e.g., private, business). \\
NS (Nights Spent) & Total number of nights spent during the trip. \\
TE (Total Expenditure) & Total monetary spending associated with the trip. \\
EC (Expenditure Categories) & Breakdown of expenses into accommodation, restaurants, and cafés. \\
\hline
\end{tabular}
\label{tab0}
\end{center}
\end{table}

\subsection{Descriptive Statistics}

To summarize the dataset, measures of central tendency and dispersion were calculated, including the mean, median, and standard deviation.

The \textit{mean} represents the average value, calculated as the sum of all values divided by the number of observations. While it provides a useful summary of the data, it is highly sensitive to outliers.

The \textit{median} is the middle value of the dataset when arranged in ascending order. Unlike the mean, the median is not affected by outliers, making it a better measure of central tendency when the data distribution is skewed.

The \textit{standard deviation} measures the dispersion of values around the mean. A higher standard deviation indicates greater variability in the dataset. It is calculated as the square root of the variance, given by Equation \ref{eq0}:

\begin{equation}{\label{eq0}}
    s = \sqrt{\frac{\sum_{i=1}^{N} (x_i - \overline{x})^2}{N-1} }
\end{equation}

where \( s \) is the standard deviation, \( x_i \) represents individual observations, \( \overline{x} \) is the mean, and \( N \) is the total number of observations.

\subsection{Visualisations and Plots}

Several graphical methods were used to explore the dataset and identify patterns.

\textbf{Scatterplots} were used to examine relationships between two continuous variables, helping to identify clusters, trends, and outliers. They are a preliminary step before calculating Pearson's correlation coefficient to ensure that the correlation reflects a true linear relationship.

\textbf{Histograms} were employed to visualise the distribution of numerical variables. They divide data into bins (intervals) and count the frequency of values within each bin. Unlike bar charts, which are used for categorical data, histograms are useful for detecting skewness (left, right, or symmetric) and assessing whether the data follows a normal distribution.

\textbf{Boxplots} (also called box-and-whisker plots) provide a summary of the data distribution. They display the minimum, first quartile (Q1), median (Q2), third quartile (Q3), and maximum values, along with potential outliers. Boxplots were particularly useful in detecting skewness, spread, and central tendency in the dataset.

\subsection{Hypothesis Testing}

To validate the assumptions of ANOVA, normality and homogeneity tests were applied using the Shapiro-Wilk and Levene’s tests.
\textbf{P-value Threshold:}  
A significance level of 0.05 was used to determine statistical significance, as defined in Equation \ref{eq1}:
\begin{equation}
p < 0.05
\label{eq1}
\end{equation}

\textbf{Shapiro-Wilk Test:}  
This test checks whether a given dataset follows a normal distribution. It is commonly used before applying parametric tests like ANOVA, t-tests, or regression analysis. A p-value greater than 0.05 indicates that the data follows a normal distribution (fail to reject the null hypothesis), whereas a p-value less than 0.05 suggests a violation of normality. While the Shapiro-Wilk test is effective for small datasets, for larger datasets, the Kolmogorov-Smirnov test is more appropriate.

\textbf{Levene’s Test:}  
Levene’s test assesses the homogeneity of variances across groups. Many statistical tests, including ANOVA and t-tests, assume that the variances of different groups being compared are roughly equal. If the p-value is greater than 0.05, homogeneity is maintained, and parametric tests can be applied. Otherwise, alternative methods such as Welch’s ANOVA or non-parametric tests are considered.

\textbf{Tukey's Honestly Significant Difference (HSD) Test:}  
A post-hoc test used after ANOVA to determine which specific group means are significantly different. It controls the Type I error rate in multiple comparisons. A p-value below 0.05 indicates a significant difference between group means.

\textbf{Fisher's Exact Test:}  
Used to assess associations between two categorical variables in a 2x2 contingency table. It is particularly useful for small sample sizes where the Chi-Square test may not be appropriate.

\textbf{ANOVA (Analysis of Variance):}  
ANOVA is used to compare the means of three or more independent groups to determine if there are significant differences between them \citep{ML_muller2024anova}. The test is applicable when data are normally distributed and homogeneity of variance is satisfied.

\textbf{Wilcoxon Test:}  
A non-parametric alternative to the paired t-test, used when data do not follow a normal distribution. The Wilcoxon-Mann-Whitney U Test is used for comparing two independent samples without assuming normality \citep{ML2_thakkar2025continuous}.

\FloatBarrier
\section{Results and Analysis}

\subsection{Tourism Data}

\subsubsection{Descriptive Statistics}

\begin{figure}[htbp]
\centerline{\includegraphics[width=0.5\textwidth]{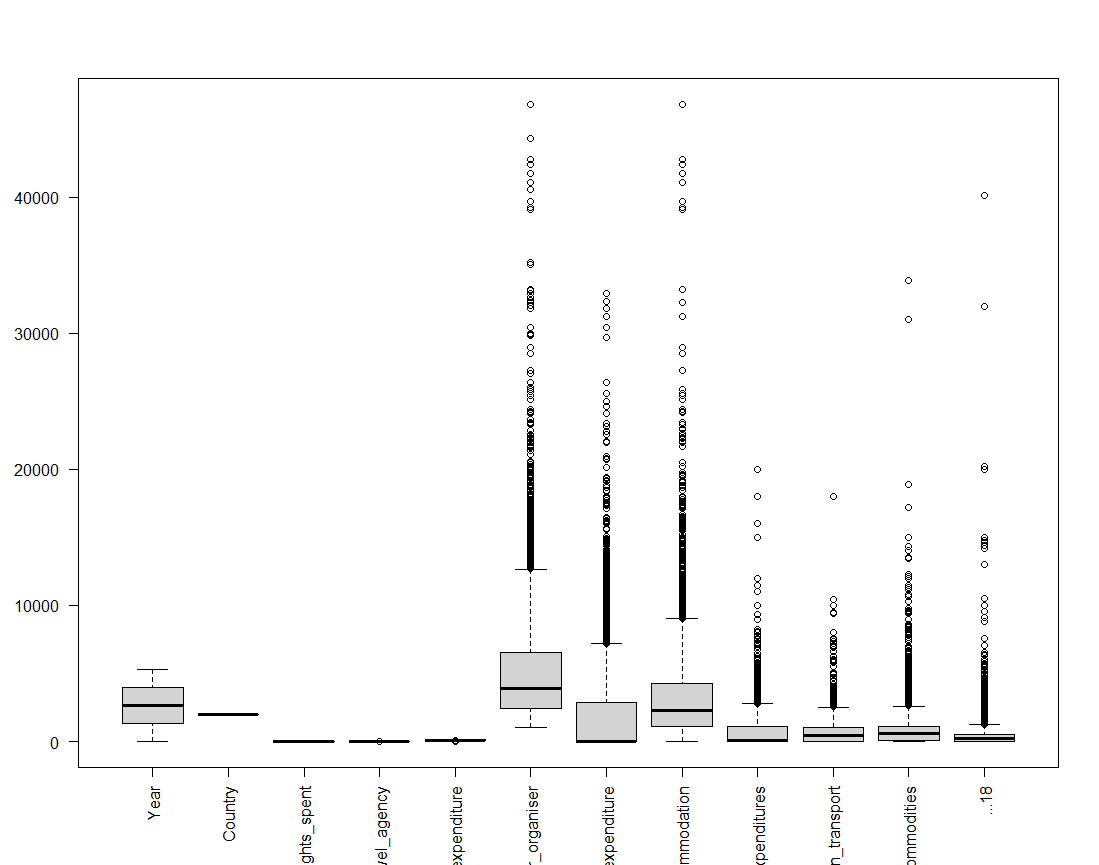}}
\caption{Box plot showing the spread of numeric variables}
\label{figbox1}
\end{figure}

A box plot (Figure \ref{figbox1}) was used to observe the variables collectively, showing the median and outliers.

A scatterplot was generated to analyse the relationship between Nights Spent and Total Expenditure. The data points cluster around specific ranges, with no clear upward or downward trend, indicating that expenditure is not directly proportional to the number of nights spent.

\begin{figure}[htbp]
\centerline{\includegraphics[width=0.5\textwidth]{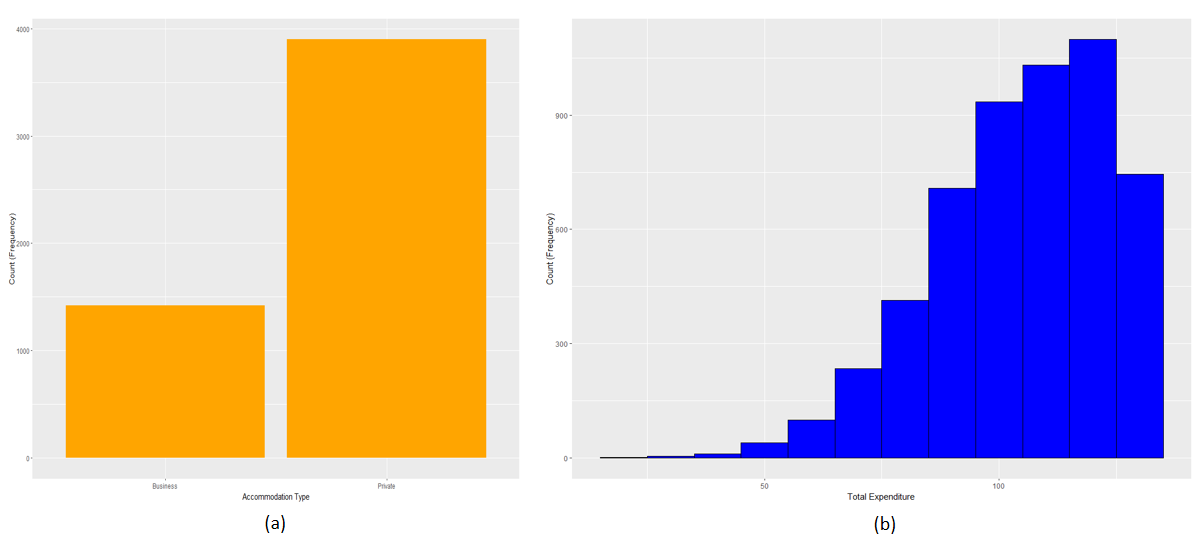}}
\caption{Histogram showing (a) distribution of Accommodation Type (b) distribution of Total Expenditure}
\label{figbox2}
\end{figure}

The histogram of Total Expenditure reveals a skewed distribution, with most spending concentrated in lower ranges, suggesting that the majority of trips had modest expenditures, with only a few high-spending outliers.

\textbf{Summary Statistics:}  
- Mean Total Expenditure: 104.86 units.  
- Median Total Expenditure: 100 units.  
- Standard Deviation: 16.8 units.  

The correlation matrix showed weak relationships among numerical variables. For example, the correlation between Nights Spent and Total Expenditure was negligible, suggesting that additional factors influence expenditure patterns.

\subsubsection{Relationships}

A correlation analysis was conducted to find relationships among the variables. Pearson’s correlation coefficient assumptions include:  
1. Both variables should be continuous.  
2. Each observation should include paired values for the two variables.  
3. Absence of extreme outliers.  
4. A linear relationship between the variables.  

If these assumptions are not met, alternative correlation measures such as Spearman’s correlation are used.

\begin{figure}[H]
\centerline{\includegraphics[width=0.5\textwidth]{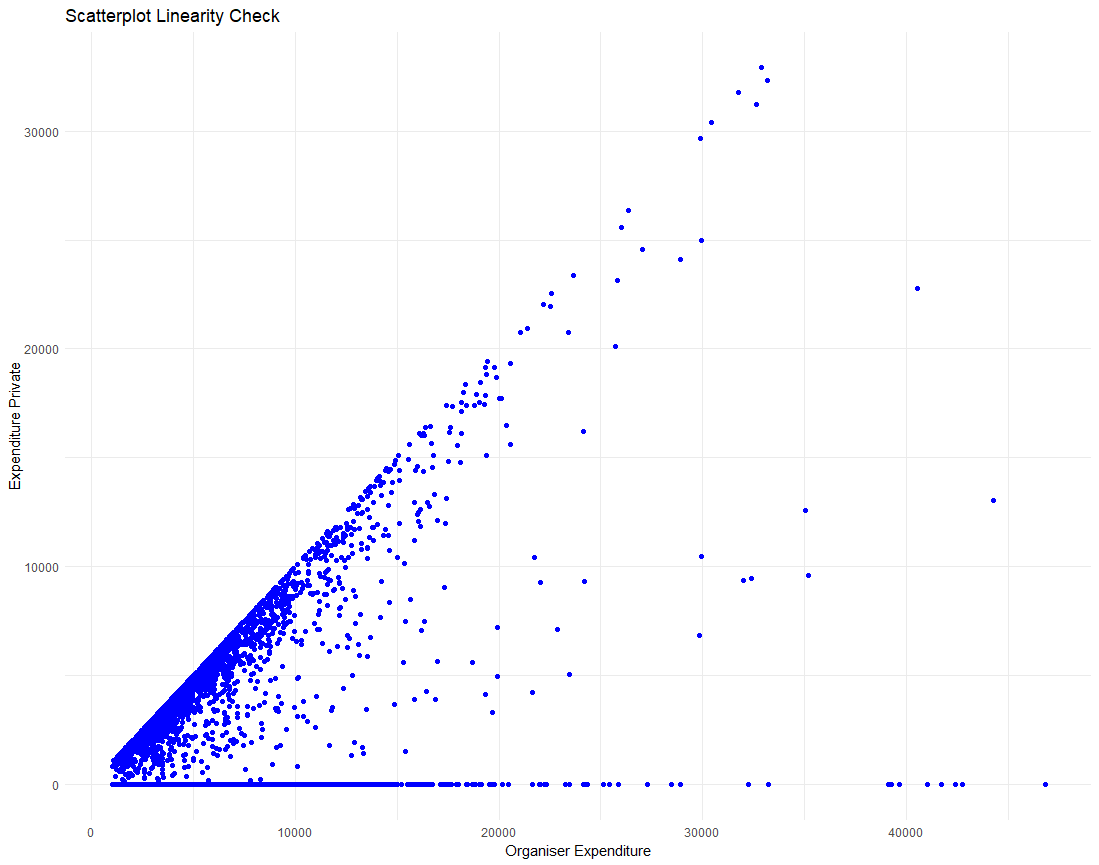}}
\caption{Scatterplot showing association between Private and Organizer Expenditure}
\label{figscatter1}
\end{figure}

A scatter plot (Figure \ref{figscatter1}) shows a strong linear relationship between Organizer Expenditure and Private Expenditure with significant anomalies.

A Pearson correlation coefficient of 0.5789 indicates a moderately strong correlation between these variables. In contrast, a weak correlation was found between Nights Spent and Total Expenditure, implying that spending is not strongly linked to the duration of stay.

\textbf{Categorical Relationships:}  
A Fisher’s Exact Test was used to assess associations between Accommodation Type and Purpose of Trip, yielding a p-value of 1e-05, indicating a significant relationship. Business travelers were more likely to stay in business accommodations, while leisure travelers preferred private accommodations.

Before conducting a T-test for total expenditure with accommodation type, homogeneity test and normality test was conducted.
Using Levene test, it was known that the groups are homogeneous. However, Shapiro-wilk test revealed that groups do not have a normality distribution. Hence, further testing would be trivial.

\section{Conclusion}

This study analysed tourism expenditure patterns 
using statistical and visualisation techniques. While significant relationships were identified among certain variables, 
ANOVA was not applicable in this case because its assumptions were not met. The key findings are summarized below:

- A significant association was found between Accommodation Type and Purpose of Trip, highlighting the influence of trip purpose on accommodation choice.

- A strong correlation was observed between Organizer Expenditure and Private Expenditure, while a weak correlation was found between Nights Spent and Total Expenditure. This indicates that spending is not strongly influenced by the duration of stay.  

- No significant difference was found in total expenditure based on accommodation type and 
trip purpose, suggesting that travelers spend similarly regardless of their trip's purpose or lodging preference.  
%
Any observed differences could be due to random variation rather than a true underlying effect.

This research provides an overview of the R programming environment and  
it contributes to the comprehension of statistical analysis methodologies and libraries in R. 
Further analysis of the dataset may reveal valuable patterns, trends, and relationships, providing deeper insights into the underlying phenomena. By applying advanced statistical and machine learning techniques, more meaningful conclusions can be derived to support data-driven decision-making.

\section*{Acknowledgment}

The author acknowledge the teaching from Dr. Dulce Gomez, J. Tiago Marques, Pedro A. Salgueiro and Sara Santos during the course: Fundamentos de análise de dados em ambiente R in Fall 2024. 
The author was supported by EURAXESS Grant ID:186817, Mobilizing Agenda: Sines Nexus, with reference "C645112083-00000059" co-financed by the PRR - Recovery and Resilience Plan.


\bibliography{refbib.bib}
\bibliographystyle{IEEEtranN}


\vspace{12pt}

\end{document}